\documentclass[a4paper,12pt]{article}
\usepackage{amsmath}
\usepackage{amssymb}
\usepackage[utf8]{inputenc}
\usepackage{graphicx}
\usepackage{subfigure}
\usepackage[font=small, skip=0pt]{caption}
\def \ni {\noindent}
\def \vs  {\vskip5mm}
\def \bea {\begin{eqnarray}}  
\def \eea {\end{eqnarray}} 
\def \mea {\nonumber\\}
\def \half {{\textstyle \frac{1}{2}}}
\def \quarter {{\textstyle \frac{1}{4}}}

\def \phip {{\Phi_+}}
\def \phim {{\Phi_-}}
\def \hatc {{\hat c}}
\def \hatE {{\widehat  E}}
\def \hatphi {{\widehat \Phi}}
\def \hatx {{\hat x}}
  
\begin{document}    
\begin{titlepage}

\title{Series solution of Painlev\'e II\\in electrodiffusion:\\conjectured convergence}     
     
\author{A.J. Bracken
and L. Bass\\
School of Mathematics and Physics\\    
The University of Queensland\\Brisbane 4072, Australia}

\date{}     
\maketitle     
\vs\ni

\begin{abstract}
A perturbation series solution is constructed with the use  of Airy functions, for a nonlinear two-point boundary-value 
problem arising in an established model of steady  
electrodiffusion in one dimension, with two ionic species carrying equal and opposite charges.  The solution 
includes a formal determination of the associated electric field, which is known to satisfy 
a form of the Painlev\'e II differential equation.  
Comparisons with
the numerical solution of the boundary-value problem   
show excellent agreement following termination of the series after a sufficient number of terms,
for a much wider range of values of the parameters in the model 
than suggested by previously presented analysis, or admitted by previously presented approximation schemes.
These surprising results 
suggest  that  
for a wide variety of cases, a convergent
series expansion for the
Painlev\'e transcendent describing the electric field  has been obtained.   A
suitable weighting of error measures for the approximations to the field and its first derivative
provides a monotonically decreasing overall measure
of the error in a subset of these cases.  It is conjectured that the series does converge for this subset.

\end{abstract}

\vs\ni

\vs\ni
{\bf E-mail:}  a.bracken@uq.edu.au, lb@maths.uq.edu.au

\end{titlepage}    

\setcounter{page}{2}

\section{Introduction}
Theoretical investigations of electrodiffusion --- the transport of charged ions in solution 
--- have a long history dating back to the pioneering works of Nernst \cite{nernst} and Planck \cite{planck}.  
Such processes underpin a variety of important 
physical and biological phenomena \cite{rubinstein1,barthel}, including 
transport of nerve impulses \cite{hodgkin1,hodgkin2,bass3}.

A well-established model in electrodiffusion 
describes one-dimensional steady transport 
of  two types of ions in solution, carrying equal and opposite charges \cite{grafov,bass1,bass2,nikonenko}. 
Charge separation of positive and negatively charged ions arising from differing diffusional fluxes, 
gives rise to an  electric field that acts back on the charges, resulting in nonlinearity of the model and associated
mathematical complexity.  
 
Grafov and Chernenko \cite{grafov} showed that in the special case when only one of the two
ionic species has a non-zero flux, across a semi-infinite domain, then
a combination of concentration and flux variables 
satisfies
 a second-order ordinary differerential equation (ODE) that is a form of Painlev\'e's 
second equation PII \cite{painleve,ince,joshi}.
Independently, one of us \cite{bass1,bass2} showed that in the general case when both fluxes can be non-zero,
another form of PII
is always satisfied by
the variable in the model that describes the  electric field, in either a semi-infinite domain or
a finite slab modeling a liquid junction.  
These independent results show  possibly the first appearances of PII 
in applications, more than $50$ years ago.  In what follows we consider further the case of the slab, 
specified by two-point boundary conditions (BCs). 

A further result obtained in 
the early studies \cite{bass1,bass2} was an approximate solution and associated asymptotic form 
for the electric field, appropriate to certain
physically well-defined limiting situations in a junction.  This result was obtained by linearizing the form of PII
satisfied by the field 
to an inhomogeneous form of Airy's ODE \cite{jeffreys,abramowitz}. 
In later work, the asymptotic behaviour of the system in other limiting situations was studied  using
the methods of singular perturbation theory \cite{macgillivray,rubinstein2,singer}.
The early results involving PII have also been extended in various other directions
\cite{zaltzman,kim}.

The form of PII found in \cite{bass1,bass2}  is quite unusual, with the feature that the coefficients 
in the ODE depend on the solution itself.  This has increased the difficulty of analytical and numerical studies, 
over and above those that apply to PII in standard form \cite{bender,clarkson1}.  
Nevertheless, existence and uniquenness
of solutions has been established \cite{thompson,park,amster1} 
in situations where physically sensible BCs are  imposed on the model.  
These ``charge-neutral" BCs require that, at each face of the slab, the concentration of the positive ions 
equals that of the negative ions.  A wide array of numerical studies has also been carried out \cite{amster1,bracken1}, 
and as a result  it can safely be concluded that under these BCs the electric field, although described by a Painlev\'e transcendent, has no singularities in the domain of interest.  

Recent research \cite{bracken1,rogers1,rogers2,amster2} has built upon the early results obtained in 
\cite{bass1,bass2}, extending known 
B\"acklund and Gambier transformations of PII \cite{gambier,lukashevich,witte}  to all the variables in the model, and exploring 
their consequences. One  development has been the discovery of the emergent phenomenon
of `B\"acklund flux quantization', wherein successive B\"acklund transformations enable the construction 
of  a sequence of exact solutions
that reflect the quantization of charge on the 
ions in solution \cite{bracken2}. 
Other studies have explored the extension of the  model to allow more ionic species in solution,
with various valences \cite{rogers3,rogers4,rogers5,amster3}.
A particular extension  to allow for fixed charges within the junction, together with
an arbitrary number of ionic species, 
has led to new differential equations in the context
of electrodiffusion,  including a generalization of PII \cite{bracken3}.

In what follows we go back to the approach initiated in \cite{bass1,bass2} and, 
building upon the linearization obtained there as 
a first approximation, 
obtain a formal solution of the whole boundary-value problem (BVP) defining the model
with charge-neutral BCs, in terms of an infinite perturbation series for the electric field,
with each term constructed using Airy functions.  The degree of approximation obtained by truncating this series
after a chosen number of terms is tested 
by comparing with the `exact' solution obtained numerically.  
These comparisons show excellent agreement for a wide range of values of 
the parameters in the model, going well beyond what is suggested by existing approximations \cite{bass1,bass2,macgillivray}.
 Most surprising is that, as the 
number of terms in the series is increased, the degree of approximation achieved
suggests that in these cases the series  converges.
If this is correct, then in this wide variety of  cases the Painlev\'e transcendent describing the steady electric field
is expressed as a convergent series, constructed using  Airy functions.    It is possible that the series provide
asymptotic, rather than convergent approximations but if so, they have a very different character, and a much wider 
range of applicability, than earlier asymptotic results \cite{bass1,bass2,macgillivray,rubinstein2,singer}.  
In any event, we see here
a development in which the power of digital computation appears to uncover unexpected possibilities for analytical approaches.

Various
connections between PII and Airy's ODE have been identified and discussed previously
\cite{ablowitz,bassom,clarkson2}.  In particular
there is a well-known sequence of exact solutions of PII expressed in terms of 
Airy functions and related by B\"acklund transformations \cite{clarkson1,rogers2}, and a link between 
this sequence and the linearization of PII obtained by Bass \cite{bass1,bass2} is known \cite{bracken1} to be 
provided with the help of the transformation first 
found by  Gambier \cite{gambier,witte}.  Nevertheless, the  relationship 
between solutions  of these two ODEs that is  suggested by our results,  is unexpected in the present context.

Another surprise is the behaviour of the errors in successive 
approximations to the Painlev\'e transcendents describing the electric field and its first derivative, 
especially in cases where the perturbation series appears to converge.  
These errors typically oscillate in size as they decrease towards zero, and 
a remarkable relationship between the oscillating 
errors in the variables describing the  field and its  derivative is revealed, which can be exploited in some cases
to obtain monotonically decreasing error measures.  We conjecture that this is a sufficient condition for the 
series in those cases to converge.  

\section{The model}
The governing system of coupled first-order ODEs for the model is 
\bea
d\,\hatc_+(\hatx)/d\hatx= ({z}e/k_B T)\,\hatE(\hatx)\,\hatc_+(\hatx)-\hatphi_+/D_+\,,
\mea\mea
d\,\hatc_-(\hatx)/d\hatx= -({z}e/k_B T)\,\hatE(\hatx)\,\hatc_-(\hatx)-\hatphi_-/D_-\,,
\mea\mea
d\,\hatE(\hatx)/d\hatx=(4\pi {z} e/\varepsilon)\left[\hatc_+(\hatx)-\hatc_-(\hatx)\right]\qquad
\label{system1}
\eea
for $0<\hatx<\delta$.  Here  $\hatE(\hatx)$ denotes
the  electric field within the junction, while $\hatc_+(\hatx)$ and $\hatc_-(\hatx)$ denote the number concentrations 
of the two ionic species;   
$\hatphi_+$ and $\hatphi_-$  their steady 
(constant) number flux-densities in the $\hatx$-direction across the junction;  $z$ the common magnitude of their 
charge numbers; 
and $D_+$, $D_-$ their diffusion coefficients. In addition,  $k_B$ denotes Boltzmann's constant,
$e$  the elementary charge, $\varepsilon$ the relative permittivity, and $T$ 
the ambient thermodynamic temperature within the solution in the junction.

An  important auxiliary quantity  is the 
electric current-density in the $\hatx$-direction, 
\bea
{\widehat J}={z}e\,(\hatphi_+-\hatphi_-)\,.
\label{current1}
\eea

As shown in Appendix A,  after all variables  are replaced by dimensionless counterparts, equations \eqref{system1}
take the form
\bea
c_+\,'(x)= E(x)\,c_+(x)-\phip\,,
\mea\mea
c_-\,'(x)= -E(x)\,c_-(x)-\phim\,,
\mea\mea
\nu E\,'(x)=c_+(x)-c_-(x)\,,
\label{system2}
\eea 
for $0< x<1$, with the prime denoting differentiation with respect to $x$, 
and \eqref{current1} becomes
\bea
J= \tau_+\phip - \tau_- \phim\,,
\label{current2}
\eea
where $\tau_+$ and $\tau_-$ are positive dimensionless constants, to be regarded as known,  with 
\bea
\tau_+ + \tau_- =1\,. 
\label{tau}
\eea

The positive dimensionless constant $\nu$ in \eqref{system2} is defined by \eqref{nu_def} below, 
and is also to be regarded as known.   
It is the squared ratio of an internal `Debye length' 
to the width  $\delta$ of the junction \cite{bass1}, and it plays an important role in the analysis of the model. 
      
Our concern in the present work is with charge-neutral boundary conditions (BCs) as in \eqref{boundary_concs}, 
supplemented by a prescribed value for the current density.  In dimensionless form, we have
\bea
c_+(0)=c_-(0)=c_0\,,\quad c_+(1)=c_-(1)=c_1 (=1-c_0)\,,\quad  J=j\,,
\label{BCs1}
\eea 
where $c_0>0$, $c_1>0$, and $j$  are prescribed constants.
Note that the system \eqref{system2} with BCs \eqref{BCs1} admits the  `reflection symmetry'
\bea
 c_{R\pm}(x) = c_{\pm}(1-x)\,,\quad\!  E_R(x) = -E(1-x)\,,\quad\!  \Phi_{R\pm} 
=- \Phi_{\pm}\,,\quad\! J_R= -J
\label{symmetry}
\eea
so there is no significant loss of generality caused by our  assumption below that 
\bea
 0<c_0<  c_1\,.
\label{conc_constraints}
\eea
(The case $c_0=c_1$ leads to a constant solution that is not of interest in what follows.)

In \eqref{system2} and \eqref{current2}, $\phip$ and $\phim$ also have constant values, 
but these values are not prescribed.  Rather, they are 
to be determined, subject to \eqref{current2} and the last of \eqref{BCs1}, as part of the solution, 
together with $c_{+}(x)$, $c_-(x)$ and $E(x)$. Thus there are
five unknowns, to be determined using the ODEs \eqref{system2} and the five BCs in \eqref{BCs1}.

\section{Painlev\'e II  boundary-value problem}
The BVP defined by \eqref{system2} and \eqref{BCs1} 
can be rewritten entirely in terms of $E(x)$ and its derivatives, following  \cite{bass1,bass2}.  
To this end, we note from \eqref{system2}  
the first-integral  
\bea
c_+(x)+c_-(x)-\half\nu E(x)^2+(\phip + \phim)x=\,\,{\rm const.}\,,
\label{first_integral}
\eea
which gives, on imposing the first two of \eqref{BCs1}, 
\bea
c_+(x)+c_-(x)=\half\nu E(x)^2
+\left\{ 2(c_1-c_0)+\half \nu [E(0)^2-E(1)^2]\right\}x
\mea\mea
+2c_0-\half\nu E(0)^2\,,\qquad\qquad\qquad\qquad\qquad\qquad\qquad
\label{first_integral2}
\eea
and
\bea
\phip + \phim=\left\{2(c_0-c_1)+\half\nu[E(1)^2-E(0)^2]\right\}\,.
\label{first_integral3}
\eea
Using the third of \eqref{system2} and the third of \eqref{BCs1}, we then have
\bea
c_{\pm}(x)= \quarter\nu E(x)^2
+\left\{ (c_1-c_0)+\quarter \nu[E(0)^2-E(1)^2]\right\}x
\mea\mea
+c_0-\quarter\nu E(0)^2 \pm \half \nu E\,'(x)\,,\qquad\qquad\qquad\qquad
\label{c_forms}
\eea
and
\bea
\Phi_{\pm}=\tau_{\mp}\left\{2(c_0-c_1)+\half\nu [E(1)^2-E(0)^2]\right\}\pm j\,,
\label{Phi_forms}
\eea
showing that all four of $c_+(x)$, $c_-(x)$, $\phip$ and $\phim$ are determined if $E(x)$ and $E\,'(x)$ can be found for all $0\le x \le 1$. 
 
To find an equation that determines $E(x)$ and $E\,'(x)$, we differentiate the third of \eqref{system2} once and 
apply the first two of \eqref{system2} and 
\eqref{first_integral} to get
\bea
\nu E\,''(x)&=&\half\nu E(x)^3  
+\left[2c_0-\half\nu E(0)^2
-(\phip + \phim)x\right]E(x)
\mea\mea
&\,&\qquad\qquad\qquad\qquad -(\phip - \phim)\,,
\mea\mea
&=& \half\nu E(x)^3 
\mea\mea
&+&\left[2c_0-\half\nu E(0)^2+\left\{2(c_1-c_0)+\half\nu [E(0)^2-E(1)^2]\right\}x\right]E(x)
\mea\mea
&+& (\tau_--\tau_+)\left\{2(c_1-c_0)+\half\nu [E(0)^2-E(1)^2]\right\}-2j\,.
\label{painleve}
\eea
From the third of \eqref{system2} and \eqref{BCs1}, we see that this second-order ODE is to be solved subject
to the Neumann BCs
\bea
E\,'(0)=0=E\,'(1)\,.
\label{BCs2}
\eea
It may be remarked that the derivation given here shows that, in effect,
the system of ODEs \eqref{system2} can be regarded as a Noumi-Yamada system \cite{noumi,matsuda} 
for \eqref{painleve}, which has been identified \cite{bass1,bass2} as a form  of PII. 
 
This apparent simplification of the BVP, now involving only 
the single second-order ODE \eqref{painleve}, comes at a cost.  The coefficients in the
ODE 
involve the unknown values $E(0)$ and $E(1)$ which have to be determined as part of its solution. 
As mentioned above, this complicates
the  theoretical and numerical analysis \cite{thompson,amster1}.  Furthermore, while
a suitable scaling of the dependent variable 
combined with
a linear transformation of the independent variable \cite{bass1,rogers1}, 
\bea
E(x)=2\beta y(z)\,,\quad z=\beta x + \gamma\,,
\label{Etrans1}
\eea
enables \eqref{painleve} to be rewritten in a standard form \cite{painleve} for PII,
\bea
y\,''(z)=2y(z)^3+zy(z)+C\,,\quad  C={\rm const.}\,,
\label{Etrans2}
\eea
the  values of the constants  $\beta$ and $\gamma$ needed to secure this result,
and consequently the value of $C$, 
all involve the unknown values 
$E(0)$ and $E(1)$.  This limits the usefulness of such a transformation and of the standard form \eqref{Etrans2} 
in the present context, and also obscures the interpretation of the results below in terms of solutions
of \eqref{Etrans2}. 

When the terms quadratic and 
cubic in the electric field are neglected in \eqref{painleve}, the equation is linearized to an inhomogeneous form of Airy's
ODE
\cite{jeffreys,abramowitz},
\bea
\nu E\,''(x)=
\left[2c_0+2(c_1-c_0)x\right]E(x)
+2\left[ (\tau_--\tau_+)(c_1-c_0)-j\right]\,.
\label{airyODE}
\eea  
It was the solution of this equation, with the BCs \eqref{BCs2}, that was used in  \cite{bass1,bass2}
to obtain the asymptotic  behaviour of the electric field in some regimes of  physical interest.  

Of particular relevance to what follows is that it was also shown there  that a sufficient condition for the validity
of this linearization is that 
\bea
\nu E_{{\rm max.}}^2 \ll 1\,,
\label{bass_condn}
\eea
where $E_{{\rm max.}}$ is an upper bound on the values  attained by $|E(x)|$ on the interval $0\leq x \leq 1$. 
 Note that $E_{{\rm max.}}$ is not
known until after the solution is (approximately) determined, so that \eqref{bass_condn} is an {\em \'a posteriori}
consistency check on the 
applicability of the linearization, rather than an indicator of when it can be applied.  

More recently, it has been shown \cite{park,bracken1} that in any solution of the BVP as defined by (3), (6),   
both $c_+(x)$ and $c_-(x)$ 
are positive for $0\leq x \leq 1$.  
Moreover,  at least one of $\phip$, $\phim$ 
is negative.  
Furthermore, it is  now known that a unique solution exists for the
reformulated BVP, and hence for the BVP as originally formulated for \eqref{system2}, 
for any given boundary values as in (6), (8) \cite{thompson,park,amster1}.

 As a result of these studies, it follows that every such solution falls into one of the following  classes:
 
\bea
&{\rm {\bf (A):}}&\; E(x)>0\,,\quad E\,'(x)<0\,,\quad c_-(x)>c_+(x)>0  \quad {\rm for}\quad 0<x<1\,,
\mea\mea
&{\rm and}&\; c_1 E(1)/c_0>E(0)>E(1)\,,\quad \Phi_-<0\,,\quad \Phi_-<\Phi_+<-\Phi_- \,; \qquad\label{sol_formsA}
\\ \mea
&{\rm {\bf (B):}}& \; E(x)<0\,,\quad E\,'(x)>0\,,\quad c_+(x)>c_-(x)>0  \quad {\rm for}\quad 0<x<1\,,
\mea\mea
&{\rm and}&\; c_1 E(1)/c_0<E(0)<E(1)\,, \quad \Phi_+<0\,,\quad  \Phi_+<\Phi_-<-\Phi_+ \,;\qquad  \label{sol_formsB}
\\ \mea
&{\rm {\bf (C):}}&\quad E(x)=0\,,\quad c_+(x)=c_-(x)=c_0+(c_1-c_0)x= c(x)\,,\,\, {\rm say},  \quad \quad\;\,
\mea\mea
&\quad&{\rm for}\quad 0\leq x\leq 1\,,\quad {\rm with}\quad \Phi_+=\Phi_- =c_0-c_1\,.\qquad\qquad\qquad\qquad
\label{exact_planck}
\eea

It can be seen from \eqref{current2} that in any such solution,  $j>0$ implies $\Phi_-<0$, and 
$j<0$ implies $\Phi_+<0$, while $j=0$ implies both
$\Phi_-<0$ and $\Phi_+<0$; but these constraints are not enough to determine the class to which 
the solution 
belongs. 

Note that it follows from \eqref{sol_formsA} and \eqref{sol_formsB} that for any solution in Class {\bf A} or {\bf B}, 
$E_{{\rm max.}}^2$ in the sufficient condition \eqref{bass_condn} is given by $E(0)^2$. 

The one solution in Class {\bf C}, which is discussed in the next section as Planck's exact solution,
is a solution of the BVP in the case when
\bea 
j=(\tau_+ - \tau_-)(c_0-c_1)=j_0\,,\quad {\rm  say}.
\label{j0_def}
\eea

\section{Series expansion}       
As noted in the Introduction, approximate solutions of the  reformulated BVP \eqref{painleve}, \eqref{BCs2},  and close variants of it, 
have been sought previously \cite{macgillivray,rubinstein2,singer} in the form of asymptotic series. That approach involves
expanding $E(x)$ in powers of  $\nu$, and 
requires the tools of singular perturbation   
theory and matched asymptotic expansions \cite{bender}, as can be seen from the fact that 
$\nu$ appears in \eqref{painleve} multiplying the highest derivative
of $E(x)$.  The zeroth order
term in the (outer) expansion so obtained 
is given by a formula for $E(x)$ first obtained by Planck \cite{planck},  in the limiting case where $\nu\to 0$ 
(the `infinite charge limit'). 

This formula forms part of Planck's approximate solution of the BVP in that limit.  Written in dimensionless form,
this approximate solution has $c_+(x)=c_-(x)=c(x)$ as in \eqref{exact_planck}, and
\bea
E(x)=\left(\Phi_+-\Phi_-\right)/2c(x)\,,\quad\Phi_{\pm}=2\tau_{\mp}(c_0-c_1)\pm j\,.
\label{planck_formula}
\eea
 
It is more important for what follows that, as can easily be checked, \eqref{planck_formula} defines an exact solution \eqref{exact_planck} of the BVP for {\em any} 
non-zero value of $\nu$, 
in the particular case when $j=j_0$ as in \eqref{j0_def}.

Our approach is to hold $\nu$ fixed at some non-zero value, and look for a series expansion in a different parameter, 
perturbing away from 
Planck's exact solution \eqref{exact_planck}. 
Thus we seek a solution of \eqref{painleve} for the case
\bea
j=j_0 +\epsilon j_1=(\tau_+ - \tau_-)(c_0-c_1)+\epsilon j_1\,, 
\label{jpert}
\eea
in the form
\bea
E(x)=0+ \epsilon\, E_1(x) + \epsilon^2 E_2(x)+\cdots\,,\qquad\qquad
\label{series1}
\eea
and write $E^{(n)}(x)$ for the series truncated after the term $\epsilon^n E_n(x)$.   
Here $\epsilon$ is a book-keeping parameter that can be set equal to $1$ after calculations are completed.  What is
relevant in the perturbation \eqref{jpert} is the value of $\epsilon j_1$.  

We note firstly from \eqref{series1} that
\bea
E(x)^2=U(x)= \epsilon^2 U_2(x)+\epsilon^3 U_3(x)+\dots
\label{Esquare1}
\eea
where $U_2(x)=E_1(x)^2$, and 
\bea
U_{2n}(x)=2E_1(x)E_{2n-1}(x)+2E_2(x)E_{2n-2}(x)+\dots  
\mea\mea
\dots + 2E_{n-1}(x)E_{n+1}(x)+E_n(x)^2\,,
\mea\mea
U_{2n-1}(x)=2\left[E_1(x)E_{2n-2}(x)+E_2(x)E_{2n-3}(x)+\dots \right.
\mea\mea
\left. \dots + E_{n-1}(x)E_{n}(x)\right]\,,
\label{Esquare3}
\eea
for $n=2\,,3\,,\dots$. Similarly
\bea
E(x)E(y)^2=V(x,y)=\epsilon^3 V_3(x,y)+\epsilon^4 V_4(x,y)+\dots\,,
\label{EEsquare1}
\eea
where
\bea
V_n(x,y)=E_1(x)U_{n-1}(y)+E_2(x)U_{n-2}(y)+\dots+E_{n-2}(x)U_2(y)
\label{EEsquare2}
\eea
for $n=3\,,4\,,\dots$. 

The ODE \eqref{painleve} can now be written with the help of \eqref{exact_planck} and \eqref{jpert} as
\bea
\nu E''(x)=2c(x)E(x)+\half\nu\left[x[V(x,0)-V(x,1)] -V(x,0)+V(x,x)\right]
\mea
\mea
-2\epsilon j_1+\half\nu(\tau_--\tau_+)[U(0)-U(1)]\,.\qquad\qquad\qquad
\label{painleve2}
\eea
We now substitute in the series expansions of $E$, $U$ and $V$, and equate terms of equal degree in $\epsilon$.  At  $\epsilon^1$ we get
\bea
\nu E_1\,''(x)- 2c(x)E_1(x)=-2j_1
= R_1\,,\,\,{\rm say}\,;
\label{E1DE}
\eea
at $\epsilon^2$
\bea
\nu E_2\,''(x)- 2c(x)E_2(x)=\half\nu (\tau_--\tau_+)[U_2(0)-U_2(1)]=R_2\,,
\label{E2DE}
\eea
and at $\epsilon^n$ for $n\geq 3$,
\bea
\nu E_n\,''(x)- 2c(x)E_n(x)=R_n(x)\,,
\label{EnDE}
\eea
where 
\bea
R_n(x)=\half\nu \left\{x[V_n(x,0)-V_n(x,1)] -V_n(x,0)+V_n(x,x)\right.
\mea\mea
\left. +(\tau_--\tau_+)[U_n(0)-U_n(1)]\right\}\,.\qquad\qquad\qquad
\label{Rn1}
\eea
The ODEs \eqref{E1DE}, \eqref{E2DE}, \eqref{EnDE} are to be  solved subject to the BCs $E_n\,'(0)=0=E_n\,'(1)$ in each case, as follows from
\eqref{BCs2}. 
 
Note that each $R_n$ depends only on the $E_k(x)$, $k<n$, determined at previous steps, so that each ODE takes the form of an Airy equation 
with known 
inhomogeneous RHS, thus permitting a solution for $E_n(x)$, $E_n\,'(x)$ 
at each stage in the form
\bea
E_n(x)={\mathcal F}_{R_n}(x)\,,\quad E_n\,'(x)={\mathcal G}_{R_n}(x)\,,
\label{Ensol}
\eea
as shown in Appendix B.  
From \eqref{c_forms} and \eqref{Phi_forms} we can then obtain series expansions for the remaining unknowns $c_{\pm}(x)$, $\Phi_{\pm}$.    

Note from
\eqref{airy7} that
\bea
\epsilon^n E_n(x)={\mathcal F}_{\epsilon^n R_n}(x)\,,\quad  \epsilon^n E_n\,'(x)={\mathcal G}_{\epsilon^n R_n}(x)\,. 
\label{Ensol2}
\eea

The approximate solution obtained in this way for the BVP of interest, up to order $\epsilon$, is 
\bea
c_+^{(1)}(x)=c_0+(c_1-c_0)x+ \epsilon \,\half \nu\, {\mathcal G}_{j_1}(x)\,,
\mea\mea
c_-^{(1)}(x)=c_0+(c_1-c_0)x- \epsilon \,\half \nu\, {\mathcal G}_{j_1}(x)\,,
\mea\mea
E^{(1)}(x)= \epsilon {\mathcal F}_{j_1}(x)\,,\qquad\qquad
\mea\mea
\Phi_+^{(1)}=c_0-c_1+\epsilon j_1\,,\quad \Phi_-^{(1)}=c_0-c_1-\epsilon j_1\,.
\label{ludvik1}
\eea
Noting \eqref{jpert} and \eqref{airy7}, we can rewrite these formulas as 
\bea
c_+^{(1)}(x)=c_0+(c_1-c_0)x+\half \nu \,{\mathcal G}_{j-j_0}(x)\,,
\mea\mea
c_-^{(1)}(x)=c_0+(c_1-c_0)x- \half \nu\, {\mathcal G}_{j-j_0}(x)\,,
\mea\mea
E^{(1)}(x)= {\mathcal F}_{j-j_0}(x)\,,\qquad\qquad
\mea\mea
\Phi_+^{(1)}=c_0-c_1- j_0+ j\,,\quad \Phi_-^{(1)}=c_0-c_1+j_0- j\,.
\label{ludvik2}
\eea
It is not hard to check that this is, in dimensionless form,  the approximate solution of the BVP \eqref{system2}, \eqref{BCs1} 
with $j_0$ as in \eqref{j0_def}, as
obtained previously \cite{bass1} 
by linearizing \eqref{painleve}, neglecting the quadratic and cubic terms on the RHS.

\section{Numerical evaluation of the series}
 {\em MATLAB} \cite{matlab} was used to conduct a variety of numerical experiments, comparing the 
 `exact' numerical solution of (3) 
with the approximate solution obtained by terminating the series in (26) at order $\epsilon^n$ 
for a variety of values $n$, for various values of the dimensionless model parameters  in the ranges 
\bea
0<\nu\leq 10\,,\quad  0<\tau_+ (=1-\tau_-)<1\,,\quad  0<c_0 (=1- c_1)<1\,,
\mea\mea
{\rm and} \quad -2.75<j (=j_0+\epsilon j_1)<2.75\,.\qquad\qquad
\label{ranges}
\eea 

It is important to emphasize at this point that in this paper we are primarily interested in the mathematical structure
of the model of binary electrodiffusion defined in Secs. 2--4,  rather than in 
the associated physics, which has been extensively discussed elsewhere in the literature 
\cite{rubinstein1,barthel,nikonenko,hodgkin1,grafov,bass1,bass2,rubinstein2,syganow}.  For this reason,
we consider $\nu$,  $\tau_+$, $c_0$ and $j$ as real 
parameters in the ranges \eqref{ranges}, 
without regard for additional constraints imposed by the physics of particular
experimental situations.  For example, while  values of $\nu$   less than $1$ have commonly been described
in the literature just cited, in particular for
nerve membranes \cite{hodgkin2,bass3}, 
a value as large as $10$ corresponds to a ratio of
internal Debye length to junction width $\delta$ of about $3.2$. 
Ratio-values greater than $1$ have also been discussed in the literature,
particularly in the context of the so-called `constant field approximation',
\cite{goldman,moore,agin,bracken3}.  They can be achieved in principle by making
the junction width 
$\delta$ or the  characteristic (reference) concentration $c_{\rm ref.}$ in (A2) sufficiently small, without
compromising charge-neutrality at junction boundaries.  Having said that, 
we recognize that
it may be difficult to achieve  ratios as large as $3.2$ in practice.    
Similar remarks apply to the values of the current-density ${\widehat J}$ 
that are implied when our dimensionless  $ |j|$  is allowed to reach values as large as $ 2.75$.  
Such  values of $\nu$ and $|j|$ together will be found to
give rise to values of $|E_{\rm max.}|$ as large as 5.5, as in  cases associated with Fig. 6 below 
for example, and so to
electric field values (in dimensional form) as large as
\bea
|\hatE_{\rm max.}|\approx 5.5  k_B T/z e \delta \,,
\label{typical_field}
\eea
leading to
$|\hatE_{\rm max.}|\delta \approx 5 \times 10^{-4}$  statVolt 
$ \approx 150$ mV for monovalent ions at room temperatures.
Depending on the value of $\delta$, field strengths of that magnitude may be very much larger than those readily 
achieved in experiments. 
    
Our justification for admitting such parameter values in the numerical exeriments that follow, 
is that we wish  to throw as much  new light as possible on the mathematical model 
of electrodiffusion and its perturbation series solution. To that end we explore a wide range of {\em mathematically sensible}
values of the parameters therein.

In what follows,  ${\mathcal S}^{{\rm num.}}$  denotes the numerical solution with components 
$c_{\pm}^{{\rm num.}}(x)$, $E^{{\rm num.}}(x)$,  {\em etc.}, 
and  ${\mathcal S}^{(n)}$ denotes  
the $n$th approximation with components  $c_{\pm}^{(n)}(x)$,   $E^{(n)}(x)$,  {\em etc.}  
Here $c_{\pm}^{(n)}(x)$ and $\Phi_{\pm}^{(n)}$  are obtained
from $E^{(n)}(x)$ and $E^{(n)}\,'(x)$ using \eqref{c_forms} and \eqref{Phi_forms}. 
(An abuse of notation is allowed below: the values of the 
components of $S^{(n)}$ appearing there are actually
approximations to the theoretical formulas of Sec. 4, obtained using numerical  evaluation of the  
integrals involved  in Appendix B.)

As remarked earlier,  the modified BVP \eqref{painleve}, \eqref{BCs2} is not straightforward to 
solve directly for $E^{{\rm num.}}(x)$
by numerical methods \cite{amster1}, 
partly because of the complication mentioned earlier,  that the coefficients in the ODE \eqref{painleve} are not 
all known {\em \'a priori}. 
This difficulty is circumvented by returning to the original BVP, rewriting it  as a system of five rather than three ODEs,
\bea
c_+\,'(x)= E(x)\,c_+(x)-\phip(x )\,,
\quad
c_-\,'(x)= -E(x)\,c_-(x)-\phim(x)\,,
\mea\mea
\nu E\,'(x)=c_+(x)-c_-(x)\,,\quad \Phi_+\,'(x)=0\,,\quad \Phi_-\,'(x)=0\qquad
\label{system3}
\eea
for $0<x<1$, with BCs
\bea
c_+(0)=c_-(0)=c_0\,,\quad c_+(1)=c_-(1)=c_1\,,\quad \tau_+ \Phi_+(0)-\tau_-\Phi_-(0)=j\,.
\label{BCs3}
\eea
This is now in a standard form, easily handled by {\em  MATLAB} \cite{matlab} library routines, 
which return numerical solutions for the unknown 
functions, including $E^{{\rm num.}}(x)$ and $E^{{\rm num.}}\,'(x)$, almost instantaneously.

In all the examples in the next Section, we take as a measure 
of the error in the $n$-th approximation
\bea
\Delta_n={\rm max}\,\left\{\left|E^{(n)}(x)- E^{{\rm num.}}(x)\right| 
+ 
\left|E^{(n)}\,'(x)-E^{{\rm num.}}\,'(x)\right|\right\}\,,
\label{Delta_def}
\eea
 bearing in mind that all components of ${\mathcal S}^{(n)}$ can be constructed from $E^{(n)}(x)$ and $E^{(n)}\,'(x)$
using \eqref{c_forms} and \eqref{Phi_forms}.
Here the maxima 
are taken over the  set of $1001$  node points  $x=0$, $0.001$, $0.002$, $\dots$, $0.999$, $1$ 
used  in determining the `exact' numerical solution ${\mathcal S}^{{\rm num.}}$. 
 
(We considered also the alternative error measure 
\bea
{\overline{\Delta}}_n&=&\left(\int_0^1\, \left\{\left[E^{(n)}(x)-E^{{\rm num.}}(x)\right]^2\right.\right.
\mea\mea
&\quad&\left.\left. + \left[E^{(n)}\,'(x)-E^{{\rm num.}}\,'(x)\right]^2\right\}\,dx\right)^{1/2}\,,
\label{Deltabar_def}
\eea
but found no substantial differences in values or behaviour between ${\overline {\Delta}}_n$ and $\Delta_n$
in the examples treated.  
Because numerical evaluation of
the integral in \eqref{Deltabar_def}  at each value of $n$  introduces further small errors,  use of $\Delta_n$ was
preferred in what follows.)

The determination of  the discretized $E^{{\rm num.}}$ and $E^{{\rm num.}}\,'$, 
and of the integrals involved in the determination of the discretized
$E^{(n)}$ and  $E^{(n)}\,'$, are all subject to numerical errors.  Adopting a conservative
position, we do not 
claim reliability of calculated values of $\Delta_n$ 
smaller than $10^{-7}$, nor do we consider such values for $n>500$.
In each example below,  we say only that the sequence
of approximations  $S^{(n)}(x)$  {\em appears to converge}
to $S^{{\rm num.}}(x)$   if a positive integer $n_7$
can be found such that $\Delta_n<10^{-7}$
for $n>n_7$, out to $n=500$, and that it {\em appears to diverge} if the size of  $\Delta_n$ 
is clearly growing in size for
sufficiently large
$n$-values, out to $n=500$.

\subsection{Examples}
In the first set of examples, $\tau_+=1-\tau_-=0.6$  and  $c_0 =1-c_1=1/3$, so that   
according to \eqref{j0_def},  

\bea
j_0= -1/15 \approx - 0.067\,.
\label{nonzero_jzero}
\eea
With these values held  fixed,  
various values of $\nu$ and $\epsilon j_1$ were considered, and in each case  the error $\Delta_n$ was 
calculated for increasing values of $n$. 

Fig. 1 shows a plot in the case $n=1$, $\nu=0.1$ and $\epsilon j_1=-0.5$ (so that $\epsilon |j_1/j_0|=7.5$).   
The `exact' numerical solution  $E^{{\rm num.}}(x)$ is shown as a solid line, and values  
of the first series approximation $E^{(1)}(x)$  
are shown at $21$ evenly spaced values of $x$ between $x=0$ and $x=1$, and marked o.

\begin{figure}[h!]
\includegraphics[trim=0.1in 3.0in 0.1in 3.5in, clip, scale=0.6]{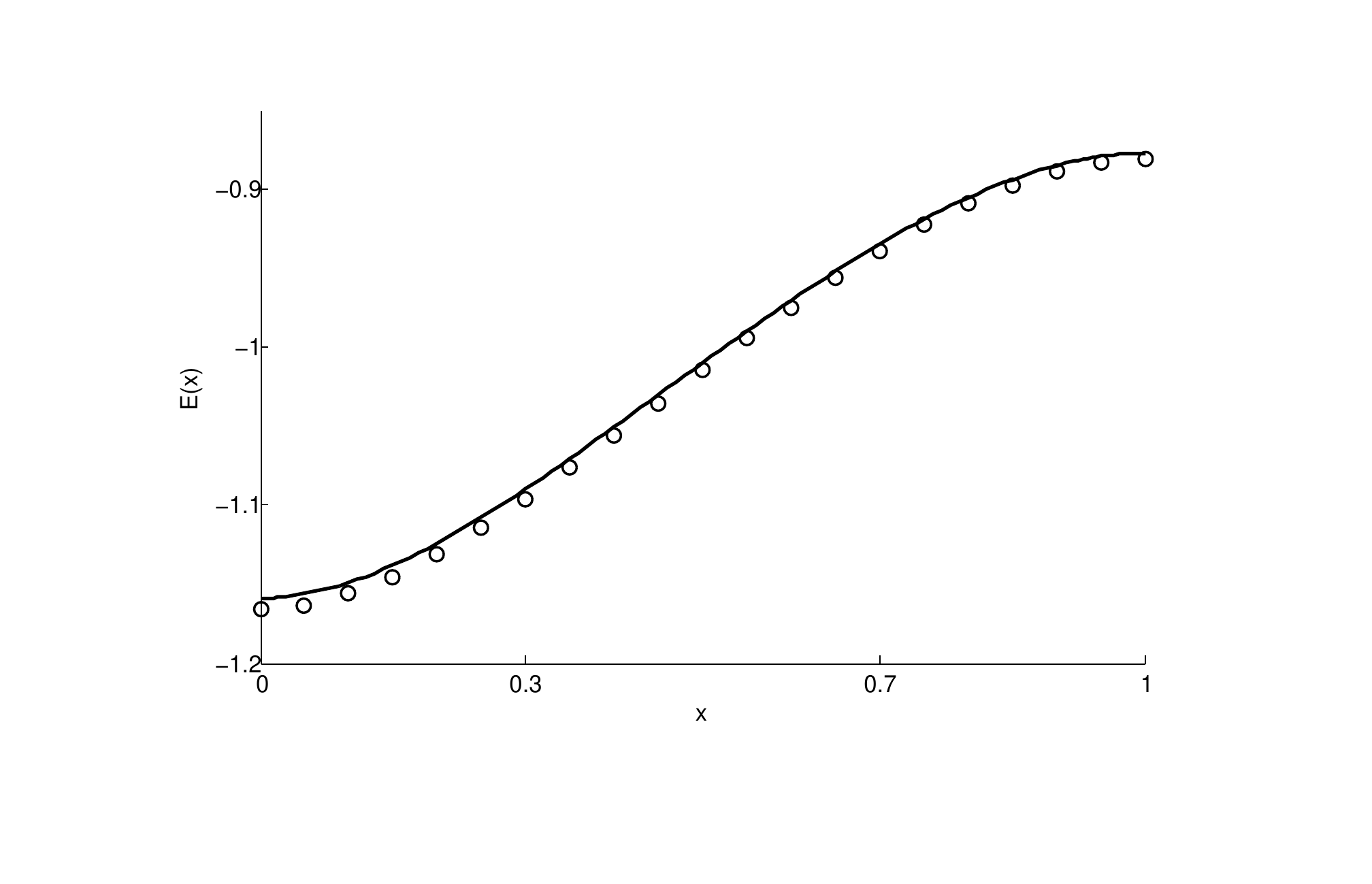}
\caption{
Plot of $E^{(1)}(x)$ (circles)
and $E^{{\rm num.}}(x)$ (solid line), showing fair agreement, in the case with  $\nu=0.1$, $\epsilon j_1=-0.5$, $\tau_+=0.6$, $c_0=1/3$.  }
\end{figure}

Fig. 2 shows, for the same settings, plots of $c_{\pm}^{{\rm num.}}(x)$, $c_{\pm}^{(1)}(x)$, 
$\Phi_{\pm}^{{\rm num.}}$ and $\Phi_{\pm}^{(1)}$.  It can be seen from the plots and \eqref{sol_formsB}  
that this case belongs to
Class {\bf B}.  
In this example, $\nu E_{{\rm max.}}^2\approx \nu E^{{\rm num.}}(0)^2\approx 0.13$ (see Fig. 1), so that the sufficient condition 
\eqref{bass_condn} is well-satisfied.  Accordingly, the  
numerical solution and this crudest ($n=1$) nontrivial approximation  \eqref{ludvik2} --- the linear approximation proposed in
\cite{bass1,bass2} --- show fair agreement, with $\Delta_1 \approx 0.02$.

\begin{figure}[h!]
\centering
\includegraphics[trim=0.1in 3.0in 0.1in 3.5in, clip, scale=0.6]{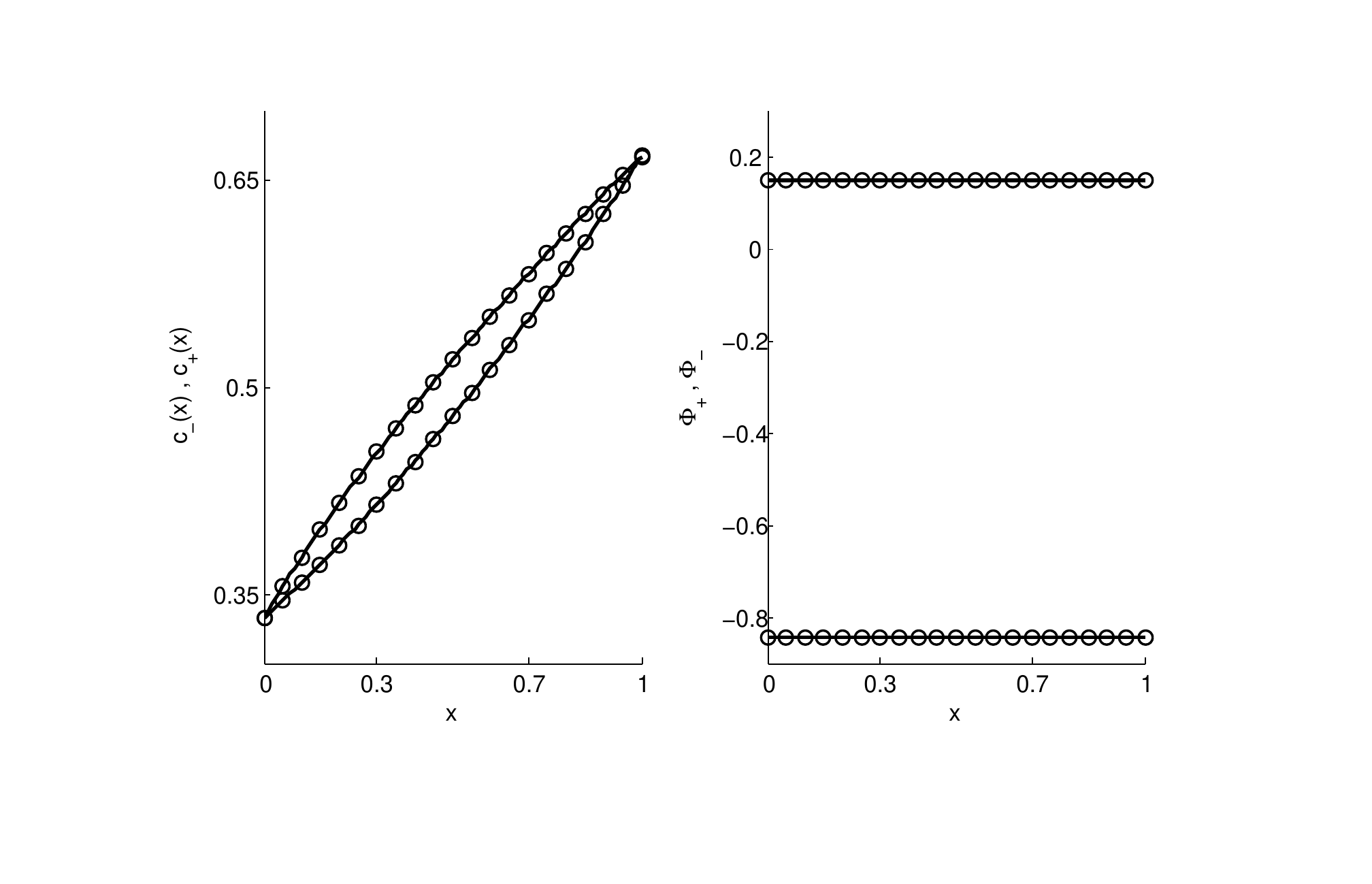}
\caption{Same constants as in Fig. 1.  On the left, plots of $c_{\pm}^{(1)}(x)$ (circles) and  $c_{\pm}^{{\rm num.}}(x)$  
(solid lines); on the right, plots of $\Phi_{\pm}^{(1)}$ (circles) and $\Phi_{\pm}^{{\rm num.}}$ (solid lines), 
showing fair agreement.}
\end{figure}

Fig. 3  on the left and in the centre shows the decrease of $\log_{10}\Delta_n$ with increasing $n$, until exceedingly 
small values 
are reached.   (Nothing should be read into the apparent constancy of the error $\Delta_n$
at values smaller than $10^{-10}$ for $n>10$, which is presumably a computational artifact.  
As remarked above, estimates smaller than $10^{-7}$ should be regarded as unreliable.  Similar remarks
apply to Figs. 4 and 5 below.)
We find 
$\Delta_n < 10^{-3}$ for $n>n_3=2$ and $\Delta_n< 10^{-7}$ for $n>n_7=7$, remaining so out to $n=500$,
and accordingly claim the the sequence of approximations  ${\mathcal S}^{(n)}$ appears to converge in this case to the solution
of the BVP defined by \eqref{system2} and \eqref{BCs1}.

Fig. 3 on the right shows plots of $E^{{\rm num.}}(x)$  and  $E^{(8)}(x)$ , 
showing the excellent agreement implied by 
 the size of $\Delta_{8}$.  Corresponding excellent agreement of 
 $c_{\pm}^{(8)}(x)$ and $\Phi_{\pm}^{(8)}$ with
 $c_{\pm}^{{\rm num.}}(x)$ and $\Phi_{\pm}^{{\rm num.}}$  follows, so we do not show the plots.

\begin{figure}[h!]
\centering
\includegraphics[trim=0.1in 3.0in 0.1in 3.5in, clip, scale=0.6]{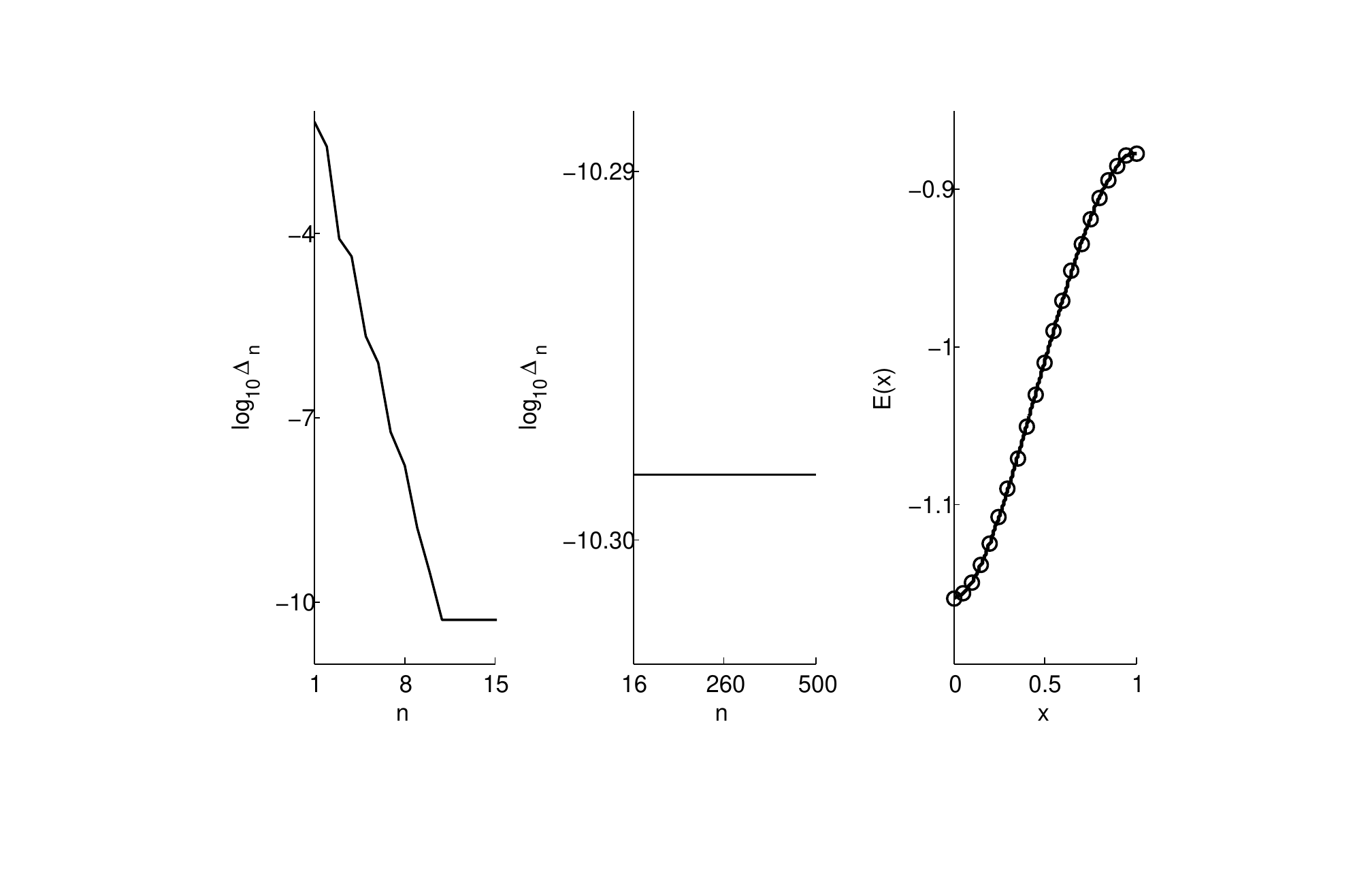}
\caption{Same constants as in Fig. 1.    On the left and in the centre, plots of 
$\log_{10}\Delta_n$ versus $n$, out to $n=500$.  On the right, plot of  $E^{(8)}(x)$ (circles)
and $E^{{\rm num.}}(x)$ (solid line), showing excellent agreement. Corresponding values of $\nu E_{\rm max.}^2$
and $\Delta_1$ are $0.13$ and $0.02$. }
\end{figure}

Lest the reader conclude that the apparent convergence here
is a consequence of the fact that   in this particular example $\nu\ll 1$  and $\nu E_{{\rm max.}}^2\ll 1$, consider next
the case where 
$\nu=1.1$ and $\epsilon j_1=-1.0$ (so that $|\epsilon j_1/j_0|=15$). 
 Fig. 4 on the left and in the centre  shows
plots of $\log_{10} \Delta_n$ versus $n$ out to $n=500$, with in this case $\Delta_1\approx 0.049$, 
$\Delta_n<10^{-3}$ for $n>n_3=4$
and $\Delta_n<10^{-7}$ for $n>n_7=11$, out to $n=500$. 
Fig. 4 on the right shows plots of  $E^{{\rm num.}}(x)$ and $E^{(12)}(x)$, 
showing the excellent agreement implied by the size of $\Delta_{12}$, and showing that this case
also belongs to Class {\bf B}. 
Not only is $\nu>1$ here, but in addition  $\nu E_{{\rm max.}}^2\approx 4.5$ (see Fig. 4), well
outside the regime defined by (19),  so the apparent convergence in this case is more surprising.

\begin{figure}[ht]
\centering
\includegraphics[trim=0.1in 3.0in 0.1in 3.5in, clip, scale=0.6]{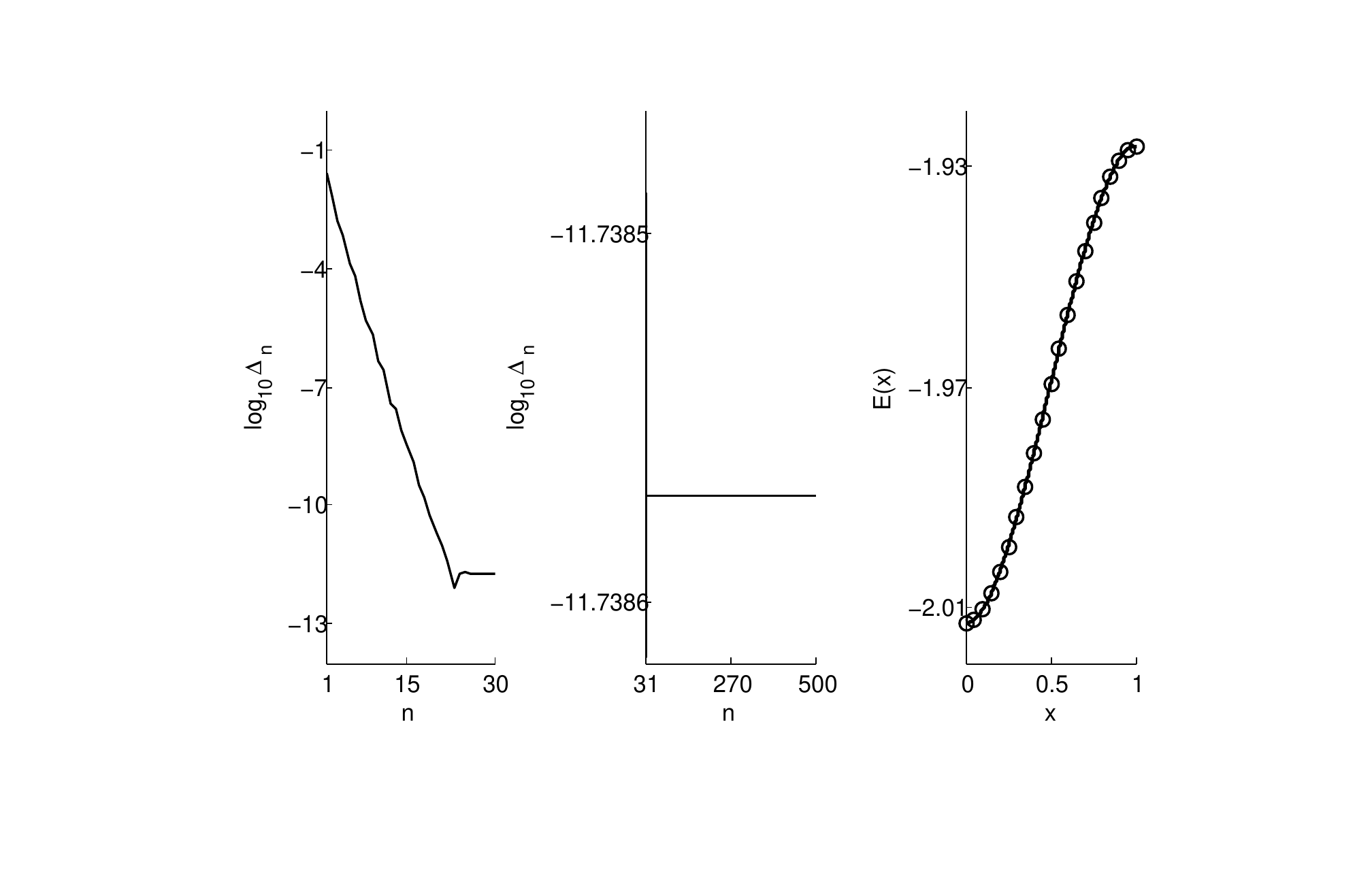}
\caption{The case with  $\nu=1.1$, $\epsilon j_1=-1.0$, $\tau_+=0.6$, $c_0=1/3$. On the left and in the centre, plots 
of  $\log_{10}\Delta_n$ versus $n$, out to $n=500$.  On the right, plot of  $E^{(12)}(x)$ (circles)
and $E^{{\rm num.}}(x)$ (solid line), showing excellent agreement. Corresponding values of $\nu E_{\rm max.}^2$
and $\Delta_1$ are $4.5$ and $0.05$. }
\end{figure}

Consider also the case where
$\nu=3.5$ and $\epsilon j_1=2.0$ (so that $|\epsilon j_1/j_0|=30$).  Fig. 5 on the left and in the centre shows
plots of $\log_{10} \Delta_n$ out to $n=500$, and in this case $\Delta_1\approx 0.17$, $\Delta_n<10^{-3}$ for $n>n_3=10$
and $\Delta_n<10^{-7}$ for $n>n_7=43$, out to $n=500$.  Fig. 5 on the right    shows plots of 
 $E^{{\rm num.}}(x)$ and $E^{(44)}(x)$ , 
showing the excellent agreement implied by the size of $\Delta_{44}$, and also showing that this case
belongs to Class {\bf A}.

 \begin{figure}[ht]
\centering
\includegraphics[trim=0.1in 3.0in 0.1in 3.5in, clip, scale=0.6]{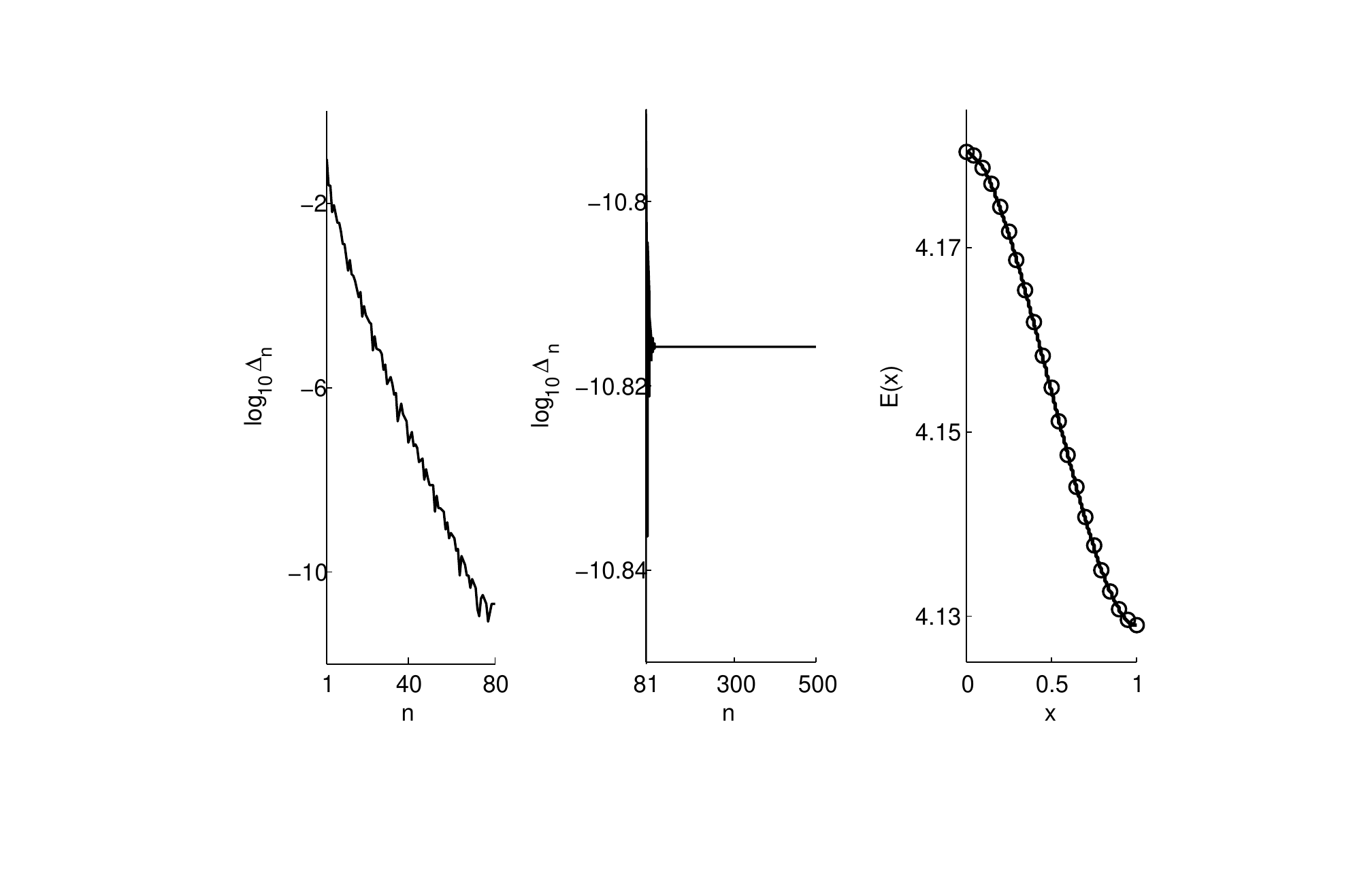}
\caption{ The case with  $\nu=3.5$, $\epsilon j_1=2.0$, $\tau_+=0.6$, $c_0=1/3$. On the left and in the centre, plots of $\log_{10}\Delta_n$ versus $n$, out to $n=500$.  On the right, plot of  $E^{(44)}(x)$ (circles)
and $E^{{\rm num.}}(x)$ (solid line), showing excellent agreement. Corresponding values of $\nu E_{\rm max.}^2$
and $\Delta_1$ are $61$ and $0.17$. }
\end{figure}

Fig. 5 also shows   that $\Delta_n$ does not decrease monotonically with
increasing $n$ in every case.  In the previous two cases, the decrease is monotonic out to $n=n_7$, 
but here that is clearly not so.  Numerous experiments show that this non-monotonic
behaviour appears to be common, rather than  exceptional.  

 Table 1 shows results from the above three cases and 
 three other apparently convergent cases.  Here as above,  $n_3$ denotes the value of
 $n$ beyond which the error $\Delta_n$ falls below $10^{-3}$, and $n_7$ the value beyond which 
 it falls below $10^{-7}$ and remains so, out to $n=500$.
 
\renewcommand{\arraystretch}{1.5}
\begin{table}[t]
\centering
\begin{tabular}[t]{|c|c|c|c|c|c|c|}
\hline
$\nu$ & $\epsilon j_1$ & $\nu E_{{\rm max.}}^2$ & $\Delta_1$ & $n_3$ & $n_7$&  Class   \\
\hline\hline  
$0.1$ & $-0.5$ & $0.13$ & $0.013$ & $ 2$ & $7$ &  $\bf{B}$\\
\hline
$0.5$ & $1.5$ & $5.2$ & $ 0.13$ & 6 &  $21$ & $\bf{A}$\\
\hline
$1.1$ & $-1.0$ & $4.5$ & $ 0.049$ & 4 &  $11$ & $\bf{B}$\\
\hline
$2.5$ & $-2.0$ & $38$ & $ 0.16$ & 11 &  $42$ & $\bf{B}$\\
\hline
$3.5$ & $2.0$ & $61$ & $ 0.17$ & 10 &  $43$ & $\bf{A}$\\
\hline
$10.0$ & $1.0$ & $42$ & $ 0.044$ & 3 &  $12$ & $\bf{A}$\\
\hline
\end{tabular}
\vs
\caption{Some apparently convergent examples with $\tau_+=0.6$, $\tau_-=0.4$, $c_0=1/3$ and $c_1=2/3$.}
\end{table} 
\renewcommand{\arraystretch}{1}

A variety of examples including those in Table 1
suggests that for the chosen values of $\tau_+$ and $c_0$, apparent
convergence occurs for any value of $\nu$ in the range $0<\nu\leq 10$ (the largest value considered) 
provided $|\epsilon j_1|<2.4$, 
but   typically becomes slower 
for a given value of $\nu$  as $|\epsilon j_1|$ 
increases towards a value in the neighbourhood of
$2.5$, and eventually fails. 

 Consider for example the behaviour of the error in the following  cases (all in Class {\bf A})  with $\nu=2.0$,
  as illustrated in Fig. 6. 
  When $\epsilon j_1 =2.45$, then
 $n_7=265$ and  the error continues to decrease (non-monotonically) out to $n=500$. 
 When $\epsilon j_1=2.48$, the situation is similar, but now $n_7=413$. 
 When $\epsilon j_1 = 2.50$, the error is still greater than $10^{-7}$ at $n=500$, though it
 is apparently still decreasing non-monotonically.  
  When $\epsilon j_1 = 2.53$, the error is still greater than $10^{-7}$ at $n=500$, and its behaviour for larger
  values of $n$ is unclear. 
 When $\epsilon j_1=2.56$,  the error decreases non-monotonically until $n\approx 110$, where it has a value
 $\Delta_{110}\approx 5\times 10^{-3}$, 
 then starts to increase non-monotonically.
 These results suggest that when $\nu=2.0$, breakdown of convergence occurs when $\epsilon j_1$ increases to a value
 somewhere between $2.48$ and $2.56$.

\begin{figure}[ht]
\centering
\includegraphics[trim=0.1in 3.0in 0.1in 3.5in, clip, scale=0.6]{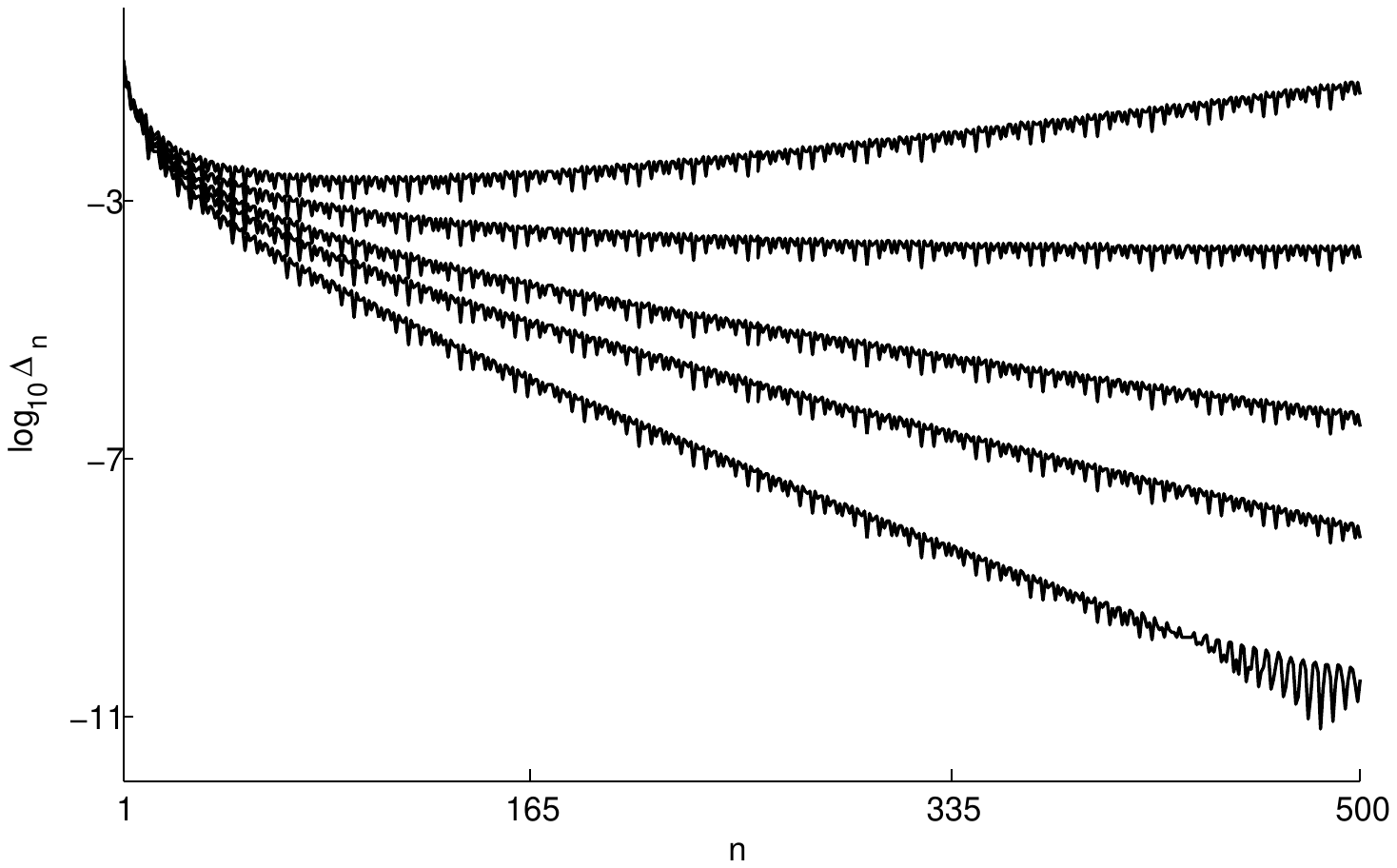}
\caption{ Plots of $\log_{10} \Delta_n$ versus $n$, out to $n=500$, 
in cases with  $\nu=2.0$, $\tau_+=0.6$ and $c_0=1/3$, and with (from bottom upwards) 
$\epsilon j_1= 2.45$, $2.48$, $2.50$, $2.53$ and $2.56$. Corresponding
values of 
$\nu E_{\rm max.}^2$ are $53$, $55$, $56$,  $57$ and $58$, with  $\Delta_1\approx 0.14$ in each case.   }
\end{figure}

The reader might suggest, after considering Fig. 6 in particular, that the series is divergent in
all of these cases --- and indeed, perhaps in every case --- and that in the apparently convergent cases,
sufficiently large values of $n$ have not been considered to show that this is so.  
Another possibility, in our view unlikely,  is that the perturbation
expansion is providing an asymptotic sequence of approximations to the exact solution in every case, 
and that this sequence gives the best
possible approximation for a particular value of $n$, say $n=n^{*}$,  with  
$n^{*}\approx 110$ in the case $\epsilon j_1=2.56$ for example,
while $n^{*}>500$ is unknown in the apparently convergent cases. 
 Our numerical experiments are unable to decide among these
possibilities, just as they are unable to prove convergence in any given case.

 Also shown in the Caption for Figs. 6  are the corresponding values of $\nu E_{\rm max.}^2$ 
 and also  of $\Delta_1$, the maximum error in the first approximation with $n=1$.  
 One of  the biggest surprises 
 in the results of this paper, strikingly illustrated by these examples, is the apparent lack of importance of the condition
$\nu E_{\rm max.}^2\ll 1$  of \eqref{bass_condn} 
 in determining whether or not the series appears to
 converge.  This is   despite the apparent importance of the size of  $\nu E_{\rm max.}^2$  in determining 
 the accuracy of the 
`linear approximation' of \cite{bass2},
 corresponding to $n=1$ here, which is reflected in the increase in the value of $\Delta_1$ with increasing
 value of the corresponding  $\nu E_{\rm max.}^2$, as is clear from the numbers shown in the Captions.    Note however that
for all the curves in Figs. 6 and 7, $\Delta_1< 1$ even though $\nu E_{\rm max.}^2\gg 1$, which is also a surprise.   
  
Consider also the behaviour in the following cases  (all in Class {\bf B}) with $\nu =1$. 
When $\epsilon j_1 =-2.45$, then
 $n_7=262$ and  the error continues to decrease (non-monotonically) out to $n=500$. 
 When $\epsilon j_1=-2.48$, the situation is similar, but now $n_7=414$. 
 When $\epsilon j_1 = -2.49$, the error is still greater than $10^{-7}$ at $n=500$, though it
 is apparently still decreasing non-monotonically.  
 When $\epsilon j_1=-2.55$,  the error decreases non-monotonically until it reaches a
 value of approximately $2.7\times 10^{-3}$ when $n\approx 135$,
 then starts to increase non-monotonically,
  much as in the previous examples illustrated in Figs. $6$.
 These results suggest in turn that when $\nu=1.0$, breakdown of convergence occurs when
 $\epsilon j_1$ decreases to a value between $-2.49$ and $-2.55$.

Examples like these suggest that the dominant influence on convergence is the value of $|\epsilon j_1|$, but
also  that the value of $\nu$
plays a minor role.  It  has not been possible to determine what, if any,   is the critical
combination of those values.

The examples above  all have  
$\tau_+=0.6$  and  $c_0 =1/3$, and the situation is more complicated
when other values of $\tau_+$ $(=1-\tau_-)$ and $c_0$  $(=1-c_1)$
are allowed, as it seems   that the values of these constants also play  roles in determining 
when the series (apparently) converges.
For example, with $\nu=1$, $\tau_+=0.9$ and $c_0 =1/3$, the series appears to converge 
when $\epsilon j_1=-2.10$ but certainly diverges when  $\epsilon j_1=-2.15$.  However, with 
 $\nu=1$, $\tau_+=0.6$ and $c_0 =0.2$,
  the series appears to converge 
when $\epsilon j_1=-2.15$ but certainly diverges when  $\epsilon j_1=-2.30$.  
For  $\nu=1$, $\tau_+=0.5$ and $c_0 =1/3$, the series appears to converge when
$\epsilon j_1=-2.5$ but certainly diverges when $\epsilon j_1=-2.75$.  Lest the reader surmise from
\eqref{jpert} that $|\epsilon j_1/j_0|$ should play the dominant role in determining convergence of the series, 
note from \eqref{j0_def} that in this last case, $j_0=0$.

In trying to determine what combination(s) of the parameters $\nu$, $\epsilon j_1$, $\tau_+ (=1-\tau_-)$ and $c_0
(=1-c_1)$ 
in the model governs (apparent) convergence,
one is faced with a ``tyranny of dimensionless constants" \cite{montroll} ---  in this case
a $4$-dimensional 
parameter-space  that has proved too large to explore fully by numerical experiments.   Similar difficuties
face efforts to determine when the error $\Delta_n$ decreases
 monotonically with increasing $n$, and when it does not.

 \subsection{Behaviour of the error}
 In order to examine more closely the peculiar oscillatory behaviour of the error
 revealed above, even in  apparently convergent cases, we consider separately now 
 the behaviour of the errors in the $n$th approximations $E^{(n)}(x)$ and $E^{(n)}\,'(x)$
 to the field and its derivative, and differently weighted sums of those errors, beyond that 
 assumed in the definition of $\Delta_n$ in \eqref{Delta_def}.    To this end, set
 
 \bea
\Delta_n(w)={\rm max}\,\left\{2w\, \left|E^{(n)}(x)-E^{{\rm num.}}(x)\right| \right.\qquad\qquad\qquad
\mea\mea
\qquad\qquad \left.+ 2(1-w)\,
\left|E^{(n)}\,'(x)- E^{{\rm num.}}\,'(x)\right|\right\}\,,
\label{Delta1_def}
\eea
where the maximum is taken over the discrete values of $x$ as in 
\eqref{Delta_def}, and
$0\leq w\leq 1$ is a weight factor at our disposal.  Note that the measure $\Delta_n$ used previously is now given by $\Delta_n(0.5)$, and that separate measures of the errors in
$E^{(n)}(x)$ and $E^{(n)}\,'(x)$ are given by $\Delta_n(1)$ and $\Delta_n(0)$, respectively. Note also that
if $\Delta_n(w)$ decreases to  values less than $10^{-7}$ out to $n=500$, for some fixed
value $0<w<1$, then we can reasonably extend our claim of apparent convergence of the 
sequences of approximations to both $E^{{\rm num.}}(x)$ and   $E^{{\rm num.}}\,'(x)$.

The values of $\Delta_n(1)$ and $\Delta_n(0)$ were computed for $1\leq n\leq n_7+1$, for each of the apparently
convergent cases in Table 1, and for the cases with $\nu=2$ and $\epsilon j_1 = 2.45$ and $2.48$ also
discussed in Sec. 5.1.  In each of these cases, it was found that 
Condition Q holds, where
\vs
\ni
{\bf Condition Q:} {\em For all} $1\leq n\leq n_7+1$,

\ni
{\em whenever} $\Delta_{n+1}(1)>\Delta_n(1)$, {\em then} $\Delta_{n+1}(0)<\Delta_n(0)$

\ni
{\em and}

  \ni  
{\em whenever} $\Delta_{n+1}(0)>\Delta_n(0)$, {\em then} $\Delta_{n+1}(1)<\Delta_n(1)$. 
\vs
This unexpected reciprocity may hold a clue to the convergence properties of the expansion.
All these cases have $\tau_+ (=1-\tau_-)=0.6$ and $c_0 (=1-c_1)=1/3$.  The condition also holds in the case
with $\nu=1$, $\tau_+=0.5$,  $c_0=1/3$ and $\epsilon j_1=-2.5$ mentioned above. 
Fig. 7 illustrates this remarkable behaviour in the case with $\nu=2.0$ and $\epsilon j_1 = 2.45$, where it
holds for $1\leq n\leq n_7 +1 =267$, as already mentioned, while the errors oscillate repeatedly in this range
during their decline
to values below $10^{-7}$.  

\begin{figure}[ht]
\centering
\includegraphics[trim=0.1in 3.0in 0.1in 3.5in, clip, scale=0.6]{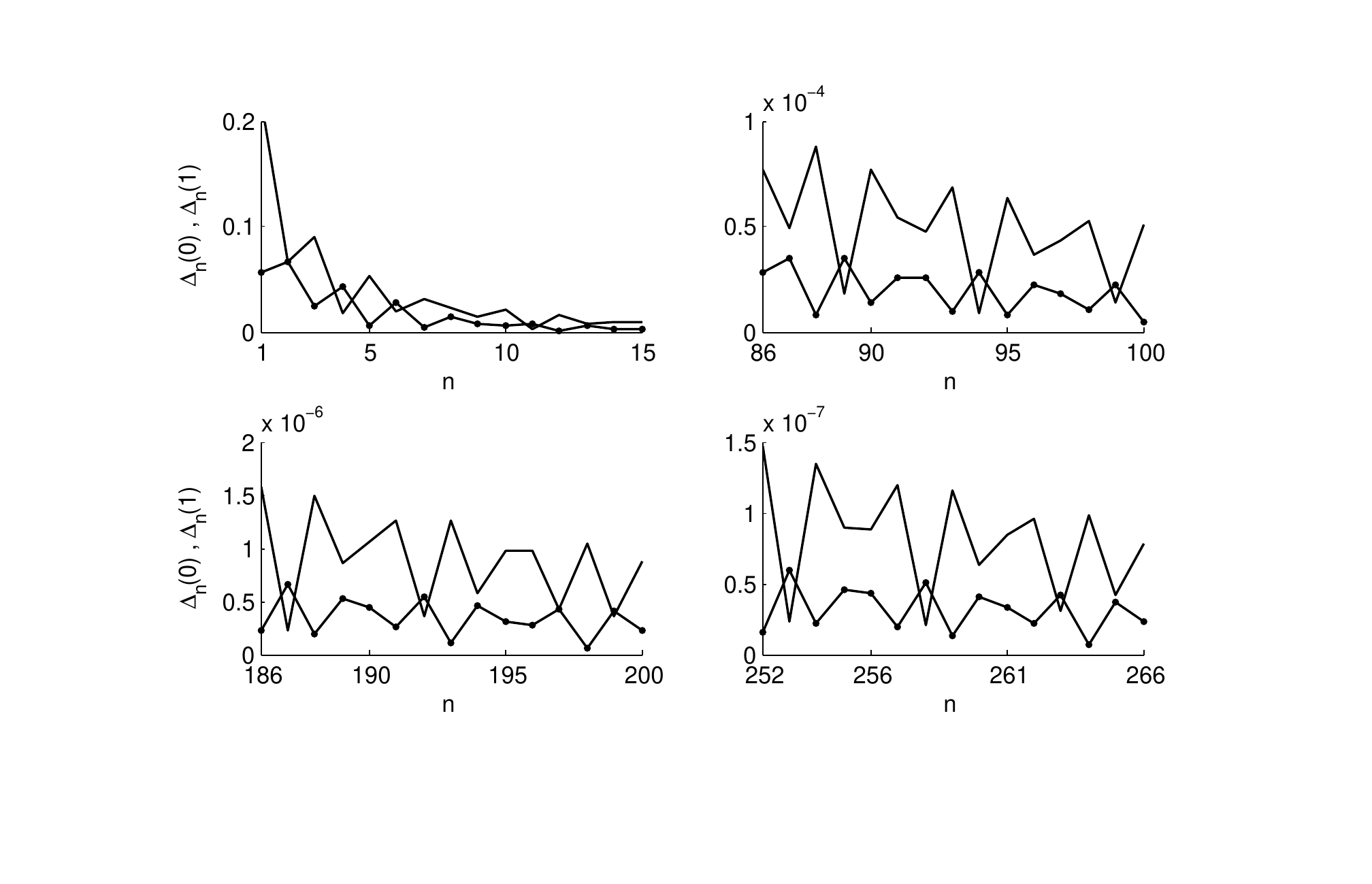}
\caption{ Plots of the errors $\Delta_n(0)$  (dot-solid) and $\Delta_n(1)$ (solid)
versus $n$ in various ranges, in the case with  $\nu=2.0$, $\tau_+=0.6$, $c_0=1/3$, and 
$\epsilon j_1=2.45$, illustrating Condition Q.
}
\end{figure}

Condition Q does not hold for all apparently convergent cases.  
It fails for the cases with $\nu=1.0$, $\tau_+=0.6$, $c_0=1/3$ and 
$\epsilon j_1=-2.45$ or $-2.48$, also discussed in Sec. 5.1, although only at a few isolated
values of $n$ in each case.  Fig. 8 shows the errors in the first of these cases for $1\leq n\leq 10$, showing 
that $\Delta_9(1)>\Delta_8(1)$ and also $\Delta_9(0)>\Delta_8(0)$.  The condition also fails to hold for the
case with $\nu=1.0$, $\tau_+=0.9$, $c_0=1/3$ and $\epsilon j_1=-2.10$,  and the case with
$\nu=1$, $\tau_+=0.6$, $c_0=0.2$ and $\epsilon j_1=-2.15$, both mentioned in Sec. 5.1 and both apparently convergent.

\begin{figure}[ht]
\centering
\includegraphics[trim=0.1in 3.0in 0.1in 3.5in, clip, scale=0.6]{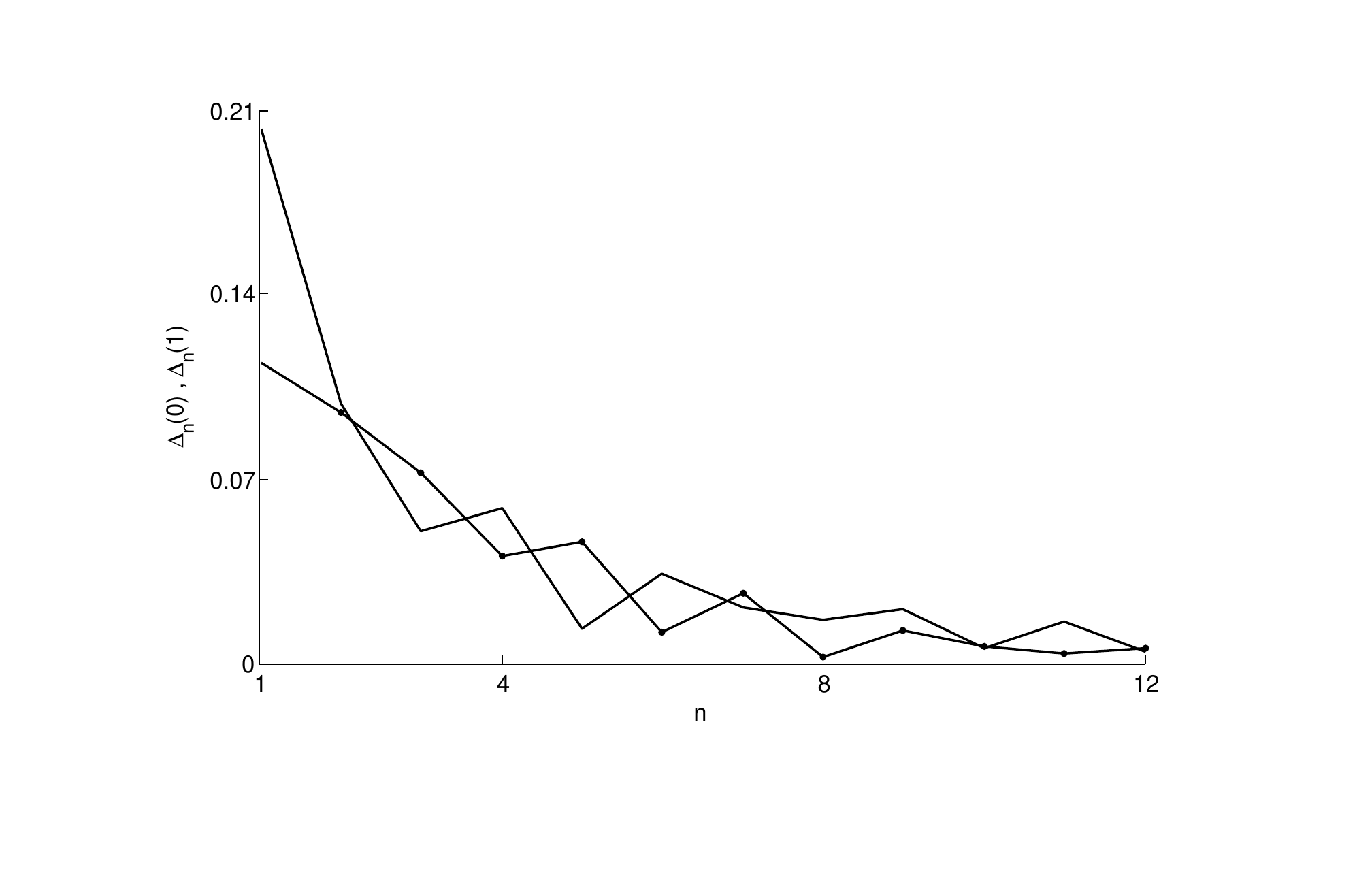}
\caption{Plots of the errors $\Delta_n(0)$  (dot-solid) and $\Delta_n(1)$ (solid) versus $n$
 in the case with  $\nu=1.0$, $\tau_+=0.6$, $c_0=1/3$, and 
$\epsilon j_1=-2.45$, illustrating the failure of Condition Q at $n=8$. }
\end{figure}

Despite these exceptions, one is led from the cases when Condition Q does hold 
to consider the possibility that a suitably weighted measure \eqref{Delta1_def}
of the 
overall error
 might in each apparently convergent case decline more smoothly towards zero than 
 \eqref{Delta_def}, perhaps even monotonically,
with the weighted up-and-down variations in  $\Delta_{n}(1)$ and $\Delta_{n}(0)$ partially cancelling
each other.  

Further numerical experiments show that a monotonic decline in the value of $\Delta_n(w)$ to less than $10^{-7}$ for 
$n>n_7$  does indeed occur in the first, third, fifth and sixth cases in Table 1, with 
$w=0.5,\,0.5,\,0.25$ and $0.2$ respectively.  In the first two of these, the weighting is just as in 
\eqref{Delta_def}, and the decline in the 
 measure of
the overall error in the field and its derivative was already seen to be monotonic in Figs. 3 and 4 .  Fig. 9 shows the graphs
of $\Delta_n(0.25)$ versus $n$ for $1\leq n\leq 15$ and $16\leq n\leq 44$ for the fifth case
in Table 1 (cf. Fig. 5), showing how the choice $w=0.25$ leads to a monotonically declining measure of the
overall error in this case.

\begin{figure}[ht]
\centering
\includegraphics[trim=0.1in 3.0in 0.1in 3.5in, clip, scale=0.6]{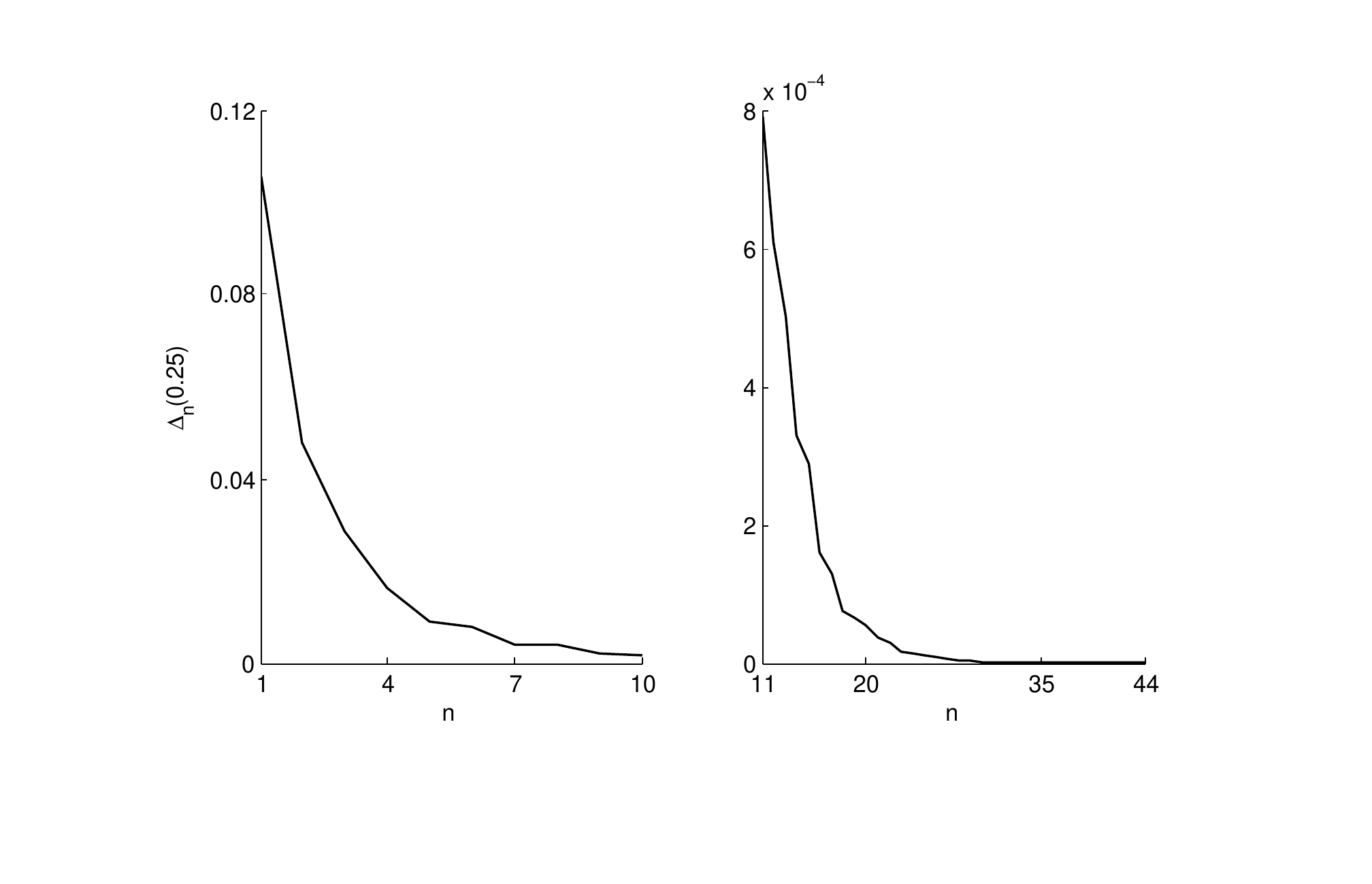}
\caption{Plots of $\Delta_n(0.25)$ versus $n$ for $n\leq 44$ in the case with $\nu =3.5$, $\tau_+=0.6$, $c_0=1/3$ and $\epsilon j_1=2.0$,  showing monotonic decline to values less than $10^{-7}$. }
\end{figure}

  It is not clear what to conclude from these numerical results.  
  We  make the
  
  \vs\ni
  {\bf Conjecture:} {\em If it is possible to  choose a weighting $0<w<1$
  in a given apparently convergent 
  case that makes the measure of error $\Delta_n(w)$ decrease monotonically to values
  less than $10^{-7}$, out to $n=500$, then the sequence of approximations $S^{(n)}(x)$ does converge
  to the exact solution of the BVP  in that case. }
  \vs
  \ni
  This leaves aside the question as to 
 what happens in other apparently convergent cases.

\section{Concluding remarks}
When the approximation to the solution of \eqref{painleve} was first obtained  \cite{bass1,bass2} by linearizing the ODE to a
form of Airy's equation --- as remarked above, this corresponds to the first
term in the formal series solution described in Sec. 4 --- the value of the result was twofold.  Firstly, it provided 
an approximate but explicit expression for the field $E(x)$, 
enabling comparisons with earlier, cruder approximations, and clarifying  
the physics of the electrodiffusion process.  Secondly, from the expression obtained for $E(x)$,
the asymptotic behaviour of the field was determined approximately at vanishingly small values 
of the current density $j$.  This was achieved 
using the known asymptotic behaviour of Airy functions \cite{abramowitz}, providing further
important information about the physics of the situation in that regime.

In these days of fast and reliable numerical computations which enable the solution of
the BVP defining the model to be obtained with great accuracy, a series expansion like
\eqref{series1} with increasingly complicated expressions for the successive terms, 
is of limited value as a means of defining better approximations than those previously obtained ---
even arbitrarily good approximations ---  to the solution of what is a complicated nonlinear
problem.   On the other hand, the expansion may have considerable theoretical significance in revealing 
an expression for the Painlev\'e transcendent satisfying \eqref{painleve} as a convergent (or asymptotic)  series, constructed
using Airy functions, for a wide range of values of the constants appearing in that equation. 

A rigorous mathematical examination of the properties of this series will be needed to determine its convergence
properties, because they
cannot be confirmed by numerical experiments 
of the type used above, however sugestive they may be.
This applies in particular to the problem of establishing the truth or otherwise of the Conjecture 
made at the end of the preceding Section.

In this context, it is worth remarking again, and may be crucial, 
that for the BCs described by \eqref{BCs1} and \eqref{conc_constraints},
the many numerical experiments that have been considered all suggest strongly that the solution
of the BVP, including the Painlev\'e transcendent electric field $E(x)$, is free from singularities
in the region $0\leq x \leq 1$ of interest.

\renewcommand{\theequation}{A\arabic{equation}}
\setcounter{equation}{0} 
\section*{Appendix A:  Dimensionless variables} 
Assuming constant boundary values for \eqref{system1}, 
\bea
\hatc_+(0)=\hatc_-(0)=\hatc_0\,,\quad \hatc_+(\delta)=\hatc_-(\delta)=\hatc_1>\hatc_0\,,
\label{boundary_concs}
\eea 
we introduce the constant reference concentration $c_{ref.}=\hatc_0+\hatc_1$ and the dimensionless constant
\bea
\nu=\epsilon k_B T/4\pi({z}e)^2\delta ^2\,c_{{\rm ref.}}\,.
\label{nu_def}
\eea.  

Then \eqref{system1} and \eqref{current1} are made dimensionless
by setting
\bea
x=\hatx/\delta\,,\quad c_{\pm}(x)=\hatc_{\pm}(\hatx)/c_{{\rm ref.}}\,,
\mea\mea
E(x)=({ z}e\delta/k_BT)\,\hatE(\hatx)\,,\qquad
\label{change_variables}
\eea
and
\bea
\Phi_{\pm}=\hatphi_{\pm}\delta /c_{{\rm ref.}} D_{\pm}\,,\quad J=\frac{\delta}{({ z}e)c_{{\rm ref.}}(D_++D_-)}\,{\widehat J}\,.
\label{change_constants}
\eea  

Setting
\bea
\tau_{\pm}= D_{\pm}/(D_++D_-)\,.
\label{tau_defs}
\eea
we obtain \eqref{system2} and \eqref{current2}.

\renewcommand{\theequation}{B\arabic{equation}}
\setcounter{equation}{0}
\section*{Appendix B:  Airy boundary-value problems} 
Consider the linear BV problem
\bea
\nu y''(x) =2c(x)\,y(x) +R(x)\,,
\quad
y'(0)=0=y'(1)\,,
\label{airy1}
\eea
where $R(x)$ is given, continuous on the interval $[0,1]$, 
and $c(x)$ is as defined in \eqref{exact_planck}.  

Linearly-independent solutions of the homogeneous ODE $(R=0)$ are provided by 
\bea
A(x)={\rm Ai}(s)\,,\quad B(x)={\rm Bi}(s)\,,\quad s=2c(x)/[4\nu(c_1-c_0)^2]^{1/3}\,,
\label{airy2}
\eea
where Ai and Bi are Airy functions of the first and second kind \cite{jeffreys,abramowitz}.  Note that
$s>0$ for all $0\leq x\leq 1$, and that Ai$(s)$ and Bi$(s)$ are both positive for $s>0$.   
Thus $A(x)$ and $B(x)$ are both positive for $0\leq x\leq 1$.  

Because the Wronskian of Ai and Bi is equal to $1/\pi$,
that of  $A$ and $B$ is given by 
\bea
W=A(x)B\,'(x)-B(x)A\,'(x)= [2(c_1-c_0)/(\pi^3\nu)]^{1/3}\,,
\label{airy3}
\eea
and  the general solution of ODE in \eqref{airy1}, according to the method of variation of parameters \cite{ince}, is then
\bea
y(x)=-\frac{1}{\nu W}\,\left\{A(x)\int_0^{x}R(y)B(y)\,dy - B(x)\int_0^{x}R(y)A(y)\,dy\right\}
\mea\mea
+ d_A\,A(x)+d_B\,B(x)\,.
\label{airy4}
\eea
with $d_A$, $d_B$ arbitrary constants.  
Imposing the BCs in  \eqref{airy1} then gives 
\bea
d_A=\frac{B'(0)}{[A'(1)B'(0)-A'(0)B'(1)]\,\nu W}\left\{A'(1)\int_0^1 R(y)B(y)\,dy \right.
\mea\mea
\left. - B'(1)\int_0^1 R(y)A(y)\,dy\right\}\,,
\mea\mea
d_B=-\frac{A'(0)}{[A'(1)B'(0)-A'(0)B'(1)]\,\nu W}\left\{A'(1)\int_0^1 R(y)B(y)\,dy \right.
\mea\mea
\left. - B'(1)\int_0^1 R(y)A(y)\,dy\right\}\,.
\label{airy5}
\eea
We write the solution $y(x)$ so defined as ${\mathcal F}_R(x)$ to emphasize its dependence on the inhomogeneous term $R(x)$. 

We note also from \eqref{airy4} that
\bea
y'(x)=- \frac{1}{\nu W}\,\left\{A'(x)\int_0^{x}R(y)B(y)\,dy - B'(x)\int_0^{x}R(y)A(y)\,dy\right\}
\mea\mea
+ d_A\,A'(x)+d_B\,B'(x)
\label{airy6}
\eea 
and, with $d_A$, $d_B$ as in \eqref{airy5}, write this function $y'(x)$ as ${\mathcal G}_R(x)$.   Finally, we note 
from the form of \eqref{airy4}, \eqref{airy5} and \eqref{airy6} that
\bea
\epsilon{\mathcal F}_R(x)={\mathcal F}_{\epsilon R}(x)\,,\quad \epsilon {\mathcal G}_R(x)={\mathcal G}_{\epsilon R}(x)\,.
\label{airy7}
\eea
for any constant $\epsilon$.

\end{document}